\documentclass[conference]{IEEEtran}
\IEEEoverridecommandlockouts
\usepackage{cite}
\usepackage{amsmath,amssymb,amsfonts,amsthm}
\usepackage{algorithmic}
\usepackage{graphicx}
\usepackage{textcomp}
\usepackage{xcolor}
\usepackage[colorlinks=true, allcolors=blue]{hyperref}
\usepackage{lineno}

\def\BibTeX{{\rm B\kern-.05em{\sc i\kern-.025em b}\kern-.08em
    T\kern-.1667em\lower.7ex\hbox{E}\kern-.125emX}}
    
\newcommand{\bmu}{\boldsymbol{\mu}}

\newcommand{\tr}{\mathrm{tr}}
\newcommand{\fH}{\mathcal{H}}

\newcommand{\toas}{\overset{a.s.}{\to}}
\newcommand{\ox}{\overline{\bx}}

\newcommand{\bX}{\mathbf{X}}
\newcommand{\bx}{\mathbf{x}}
\newcommand{\bd}{\mathbf{d}}
\newcommand{\bfr}{\mathbf{R}}
\newcommand{\bfrh}{\hat{\mathbf{R}}}
\newcommand{\bU}{\mathbf{U}}
\newcommand{\bu}{\mathbf{u}}

\newcommand{\bS}{\mathbf{S}}

\newcommand{\bw}{\mathbf{w}}

\newcommand{\tT}{T_{\mathrm{LW}}}
\newcommand{\oZ}{Z}
\newcommand{\bI}{\mathbf{I}}
\newcommand{\bfrhlw}{\bfrh_{\mathrm{LW}}}
\newcommand{\bD}{\mathbf{D}}
\newcommand{\myS}{\epsilon}
\newcommand{\mya}{\alpha}
\newcommand{\myb}{\boldsymbol{\beta}}

\theoremstyle{definition}

\newcommand{\bsalg}{BS96}
\newcommand{\cqalg}{CQ10}
\newcommand{\lalg}{LAPPW20}

\begin{document}

\title{An Improvement on the Hotelling $T^2$ Test Using the Ledoit-Wolf Nonlinear Shrinkage Estimator\\
}

\makeatletter
\newcommand{\linebreakand}{%
  \end{@IEEEauthorhalign}
  \hfill\mbox{}\par
  \mbox{}\hfill\begin{@IEEEauthorhalign}
}
\makeatother

\author{\IEEEauthorblockN{Benjamin Robinson}
\IEEEauthorblockA{\textit{Sensors Directorate} \\
\textit{US Air Force Research Laboratory}\\
Dayton, Ohio \\
ORCID: 
0000-0002-9391-4543}
\and
\IEEEauthorblockN{Robert Malinas}
\IEEEauthorblockA{\textit{EECS} \\
\textit{University of Michigan}\\
Ann Arbor, MI, USA \\
ORCID: 0000-0002-1318-5564}
\linebreakand
\IEEEauthorblockN{Van Latimer}
\IEEEauthorblockA{\textit{Mathematics} \\
\textit{UCLA}\\
Los Angeles, California, USA \\
ORCID: 0000-0001-9789-7739} 
\and
\IEEEauthorblockN{Beth Bjorkman Morrison}
\IEEEauthorblockA{\textit{Sensors Directorate} \\
\textit{US Air Force Research Laboratory}\\
Dayton, Ohio \\
ORCID: 0000-0001-7406-0311}
\and
\IEEEauthorblockN{Alfred O. Hero, III}
\IEEEauthorblockA{\textit{EECS} \\
\textit{University of Michigan}\\
Ann Arbor, MI, USA \\
ORCID: 0000-0002-2531-9670}
}

\maketitle

\begin{abstract}
Hotelling's $T^2$ test is a classical approach for discriminating the means of two multivariate normal samples that share a population covariance matrix. Hotelling's test is not ideal for high-dimensional samples because the eigenvalues of the estimated sample covariance matrix are inconsistent estimators for their population counterparts. 
We replace the sample covariance matrix with the nonlinear shrinkage estimator of Ledoit and Wolf 2020. We observe empirically for sub-Gaussian data that the resulting algorithm dominates past methods (Bai and Saranadasa 1996, Chen and Qin 2010, and Li et al. 2020) for a family of population covariance matrices that includes matrices with high or low condition number and many or few nontrivial---i.e., spiked---eigenvalues.

\end{abstract}


\section{Introduction} \label{sec:introduction}

A fundamental problem in statistics and signal processing is determining whether two independent samples have the same mean. For multivariate samples that are Gaussian and have a shared population covariance matrix, classical methods like the Hotelling $T^2$ test apply \cite{anderson1963asymptotic}. However, there is no standard technique if the sample dimension is substantial compared to the sample sizes: the so-called large-dimensional regime.

Hotelling's $T^2$ test relies on estimating the shared population covariance matrix using the sample covariance matrix. In the large-dimensional regime, the sample covariance matrix's eigenvalues are  inconsistent estimators for their population counterparts \cite{paul2007asymptotics}, leading to poor performance of Hotelling's test.
Additionally, the sample covariance matrix can be considerably more ill-conditioned than the population covariance matrix, resulting in numerical instabilities \cite{johnstone2001distribution}.  
As a result, several authors have proposed alternatives to the standard $T^2$ test.  Bai and Saranadasa \cite{bai1996effect} and Chen and Qin \cite{chen2010two} propose tests (\bsalg{} and \cqalg{}) that are an improvement for a well-conditioned population covariance matrix, and Li et al. \cite{li2020adaptable} propose a test (\lalg{}) that is an improvement if the population covariance follows the spiked model of Johnstone \cite{johnstone2001distribution}. To the best of our knowledge, until now no test improves upon of all these methods under more general assumptions.



In this paper, we propose a replacement for Hotelling's $T^2$ test and present some simulations in which it dominates \bsalg{}, \cqalg{}, and \lalg{}. We do not restrict the condition number or the number of spiked eigenvalues of the population covariance matrix in our testing.
Our method, similar to \lalg{}, is to replace the sample covariance matrix in Hotelling's $T^2$ with a matrix that has a smaller condition number. More precisely, we use the covariance matrix estimator of Ledoit and Wolf \cite{ledoit2020analytical}, which belongs to Stein's shrinkage class \cite{stein1975estimation,stein1986lectures}. By contrast, \lalg{} uses a diagonal-loading estimator, which also belongs to Stein's class, but is simply a sum of the sample covariance matrix and a scalar multiple of the identity.  
We argue for our test's optimality within Stein's class using ideas similar to those some of us have applied to other detection problems \cite{robinson2021space}. 
We note that we are not the first to have thought of applying Ledoit-Wolf-type estimators to two-sample testing, but are developing the idea toward maturity \cite{namdari2021high}.



In Section~\ref{sec:background}, we explore past improvements on Hotelling's $T^2$ test. In Section~\ref{sec:lw}, we give the Ledoit-Wolf estimator utilized in our approach, present the proposed test, and argue for the test's asymptotic optimality and predict a precise asymptotically constant false-alarm rate. In Section~\ref{sec:sim}, we provide simulations that show empirical improvement over past methods. Finally, in Section~\ref{sec:conclusion}, we present our conclusions.

\section{Background}\label{sec:background}


Suppose we have independent random $p$-dimensional column vectors $\bx_{ij}\sim \mathcal{N}(\bmu_i, \bfr)$ for $j=1,2,\dots, n_i$ and  $i\in\{1,2\}$, where $\bmu_i$ are unknown means and $\bfr$ is an unknown $p\times p$ symmetric positive-definite population covariance matrix.
We want to test  the following hypotheses
\begin{equation*}
    \begin{cases}
    \fH_0: & \bmu_1 = \bmu_2, \\
    \fH_1: & \bmu_1 \ne \bmu_2.
    \end{cases}
\end{equation*}

If the covariance matrix $\bfr$ is known, a reasonable detector is the Mahalanobis-distance detector
\begin{equation} \label{eq:quadratic-detector}
     (\ox_1-\ox_2)'\bfr^{-1}(\ox_1-\ox_2) \underset{\fH_0}{\overset{\fH_1}{\gtrless}} \tau,
\end{equation}
where $(\, \cdot\, )'$ denotes the transpose, $\tau \in (0,\infty)$, and $\ox_i$ is the sample mean of the $\bx_{ij}$'s. For unknown $\bfr$, the classical replacement of \eqref{eq:quadratic-detector} is Hotelling's $T^2$ test:
\begin{equation} \label{eq:pooled-scm}
    (\ox_1-\ox_2)'\bS_n^{-1}(\ox_1 - \ox_2) \underset{\fH_0}{\overset{\fH_1}{\gtrless}} \tau,
\end{equation}
where $n=n_1+n_2-2$ and $\bS_n$ is the ``pooled'' sample covariance matrix, given by
\begin{equation} \label{eq:scm}
\bS_n = \frac{1}{n}\sum_{i=1}^2\sum_{j=1}^{n_i} (\bx_{ij}-\ox_i)(\bx_{ij}-\ox_i)'.
\end{equation}
 The performance of the above test is well-characterized
when $n_1, n_2 \to \infty$ and $p$ is fixed \cite{anderson1963asymptotic}; however, this test can become inadmissible in the regime where $p, n_1, n_2 \to\infty$ and $p/n_1 \rightarrow \gamma_1$ and $p/n_2 \rightarrow \gamma_2$ for some $\gamma_1, \gamma_2 \in (0,\infty)$ \cite{bai1996effect}. This latter limit is known as the \emph{high-dimensional asymptotic regime}.


Bai and Saranadasa \cite{bai1996effect} suggest an improvement, \bsalg{}, under the assumption that $\left\Vert \bfr \right\Vert = o\left(\sqrt{\tr\left(\bfr^2\right)}\right)$.  This assumption yields that the matrix $\bfr$ is well-conditioned, for otherwise the two referenced quantities would be similar in size.  Their proposed test statistic is
\begin{equation}
    \frac{\frac{n_1 n_2}{n_1+n_2}\left\Vert\ox_1-\ox_2\right\Vert^2 - \tr \bS_n}{\sqrt{\frac{2(n+1)}{n}}B_n},
\end{equation}
where
\[
B_n = \frac{n^2}{(n+2)(n-1)}\left(\tr\left(\bS_n^2\right)-\frac{1}{n}\left(\tr \bS_n\right)^2\right).
\]
The corresponding detector is desirable in that, asymptotically, it has known constant false-alarm rate and predictable detection power.

In response to the need to analyze ultra-high dimensional data, Chen and Qin \cite{chen2010two} propose a test, \cqalg{}, that is similar to \bsalg{} but is additionally well suited even to values of $n$ much smaller than $p$.  Their test possesses desirable provable optimality properties, provided again an assumption that implies a low condition number: $\tr(\bfr^4) = o\{\tr^2(\bfr^2)\}$.  Their test statistic is
\[
\frac{\sum_{i\ne j}^{n_1} \bx_{1i}'\bx_{1j}}{n_1(n_1-1)} + \frac{\sum_{i\ne j}^{n_2} \bx_{2i}'\bx_{2j}}{n_2(n_2-1)} - \frac{2\sum_{i=1}^{n_1}\sum_{j=1}^{n_2} \bx_{1i}'\bx_{2j} }{n_1 n_2}.
\]

Li et al. \cite{li2020adaptable} assume instead that $\bfr$ follows the well-known \emph{spiked covariance model} of Johnstone \cite{johnstone2001distribution}, i.e., that all the eigenvalues of $\bfr$ are equal to one except for a fixed finite number that are larger than one. This assumption potentially contrasts with the well-conditionedness assumption of \bsalg{} and \cqalg{}.
Li et al. propose a test, \lalg{}, that replaces $\bS_n$ in \eqref{eq:pooled-scm} by the \emph{diagonal-loading} estimator $\bS_n+\lambda\bI$, for some optimal \emph{loading factor} $\lambda > 0$.  Like \bsalg{} and \cqalg{}, \lalg{} has a known asymptotically constant false-alarm rate and predictable detection power. Further, Li et al. remark that 
\lalg{} reduces to \bsalg{} for population covariance matrices that are both well-conditioned and spiked.

The covariance estimator of \lalg{} belongs to the shrinkage class of Stein \cite{stein1975estimation,stein1986lectures}.
These are estimators that differ from $\bS_n$ only in their eigenvalues. That is, all of the eigenvectors of the estimator are eigenvectors of $\bS_n$, but the eigenvalues may be different from $\bS_n$'s. Shrinkage covariance estimators have been studied in the high-dimensional regime in the spiked model by Donoho et al. \cite{donoho2018optimal}, and in a more general model by Ledoit and Wolf  \cite{ledoit2011eigenvectors,ledoit2012nonlinear,ledoit2018optimal}.  In particular, Ledoit and Wolf have devised closed-form expressions for shrinkage eigenvalues that are asymptotically optimal with respect to many criteria, including Stein's loss, inverse Stein's loss, Frobenius loss, inverse Frobenius loss, and so-called minimum-variance loss \cite{ledoit2020analytical}.
In what follows, we will propose again replacing $\bS_n$ in Hotelling's $T^2$ test, but rather than a diagonal-loading estimator, we utilize a shrinkage estimator of Ledoit and Wolf. 

\section{Proposed Test} \label{sec:lw}

In this section, we define our proposed test, give an indication of its asymptotic false-alarm rate, and argue for a type of asymptotic detection-theoretic optimality.

\subsection{Definition} \label{sec:def}

In order to define our proposed test, we must first describe Ledoit and Wolf's nonlinear shrinkage estimator from \cite{ledoit2020analytical}.  Throughout this paper, we will consider the high-dimensional limit in which $n,p\to\infty$ and $p/n=p/(n_1+n_2-2)\to\gamma\in (0,\infty)$. The notation $\to$ will always refer to the limit in which $n,p\to\infty$ and $p/n\to\gamma$. For notational convenience, we will follow the convention of \cite{benaych2016lectures} that quantities that are not explicitly deemed constant are varying with $n$ and $p$.  Identifying $\bS_n$ with $\bS$, then, we let $\bS = \bU\mathrm{diag}(\boldsymbol{\lambda})\bU'$ be an eigen-decomposition, where for emphasis $\bU$ and $\boldsymbol{\lambda}$ depend on $n$ and $p$.
We list the eigenvalues $\boldsymbol{\lambda}=(\lambda_1,\dots,\lambda_p)$ in non-increasing order, with corresponding eigenvectors $\bU=[\bu_1, \bu_2, \dots, \bu_p]$.

The Ledoit-Wolf nonlinear shrinkage estimator can be described as follows. Let $h_j = n^{-1/3}\lambda_j$ and $[y]^+ = \max\{y, 0\}$.  Define
\begin{align*}
a(\lambda, \boldsymbol{\lambda}) & := \sum_{j=[p-n]^+ +1}^{p}\left\{ -\frac{3(\lambda - \lambda_{j})}{10\pi h_{j}^2} +  \right. 
 \left. \frac{3}{4\sqrt{5}\pi h_{j}}\times \right. \\
 &\qquad \left. \left[1 - \frac{1}{5}\left( \frac{\lambda-\lambda_{j}}{h_{j}} \right)^2 \right] 
 \log\left| \frac{\sqrt{5} h_{j} - \lambda + \lambda_{j}}{\sqrt{5} h_{j} + \lambda - \lambda_{j}} \right| \right\},
\end{align*}
and
\[
b(\lambda, \boldsymbol{\lambda}) := \sum_{j=[p-n]^+ +1}^{p} \frac{3}{4\sqrt{5}h_{j}} \left[ 1- \frac{1}{5}\left( \frac{\lambda - \lambda_{j}}{h_{j}} \right)^2\right]^+,
\]
and
\[
s(\lambda, \boldsymbol{\lambda}) = \pi(a(\lambda,\boldsymbol{\lambda})+ib(\lambda,\boldsymbol{\lambda}))/\min\{n,p\}.
\]
Next, the shrunken eigenvalues $\hat{d}_i$ are defined as
\begin{equation*}
	\hat{d}_i := \begin{cases}
	\frac{\lambda_{i}}{\left|1-p/n-(p/n)\lambda_{i}s(\lambda_i, \boldsymbol{\lambda})\right|^{2}}, & \lambda_{i}>0\\
	\frac{1}{(p/n-1)a(0,\boldsymbol{\lambda})/n}, & \lambda_{i}=0.
	\end{cases}\label{eq:shrinkage-formula}
	\end{equation*}
(We assume $p\ne n$.) Finally, the Ledoit-Wolf estimator, denoted $\bfrhlw$, is given by $\bfrhlw = \bU\mathrm{diag}(\hat{d}_1, \hat{d}_2, \dots, \hat{d}_p)\bU'$.  

We define our test to be the Hotelling $T^2$ test with $\bfrhlw$ substituted for $\bS$:
\begin{equation*} 
\tT^2 := \frac{n_1 n_2}{n_1+n_2}(\ox_1-\ox_2)'\bfrhlw^{-1}(\ox_1-\ox_2) \underset{\fH_0}{\overset{\fH_1}{\gtrless}} \tau.
\end{equation*}
For the false-alarm rate analysis that follows, we shift and re-scale $\tT^2$ for the equivalent test
\begin{equation}
    \label{eq:oZ}
\oZ = \frac{1}{\sqrt{2p}}(\tT^2 - p) \underset{\fH_0}{\overset{\fH_1}{\gtrless}} \tau.
\end{equation}

\subsection{False-Alarm Rate}

A key consideration about any test statistic is whether its null distribution can be characterized.
In the following, we argue intuitively that the null distribution of $Z$ should be asymptotically standard normal for Gaussian data and provide empirical support for this assertion.  We expect a similar result holds for sub-Gaussian data.

Let $\sigma_i^2 = \bu_i'\bfr\bu_i$. The test in subsection~\ref{sec:def} is related to the test statistic $\tilde{Z}$, defined as follows.  Let $\bfrh$ be a shrinkage estimator $\bU\bD\bU'$, where without loss of generality, 
\[
\bD = \mathrm{diag}\left(\frac{\sigma_1^2}{c_1}, \frac{\sigma_2^2}{c_2}, \dots, \frac{\sigma_p^2}{c_p}\right),
\]
for some positive coefficients $\{c_i$\}.  Define
\begin{equation} \label{eq:oZ0}
\tilde{Z} := \frac{1}{\sqrt{2\sum_i c_i^2}}\left(\tilde{T}^2 - \sum_{i=1}^p c_i\right),
\end{equation}
where $\tilde{T}$ is given by
\begin{equation} \label{eq:ttilde0}
\tilde{T}^2 := \frac{n_1 n_2}{n_1 + n_2}(\ox_1-\ox_2)'\bfrh^{-1}(\ox_1-\ox_2).
\end{equation}
By the Berry-Esseen theorem, the conditional distribution of $\tilde{Z}$ given $\bS$ is almost surely asymptotically standard normal.

The shrunken eigenvalues $\hat{d}_i$ of Ledoit and Wolf are designed to approximate $\sigma_i^2$.  A weak form of this approximation that holds under a general matrix model laid out \cite{ledoit2011eigenvectors,ledoit2020analytical} is that for all fixed intervals $I \subset [0, \infty)$
\begin{equation} \label{eq:lwp-law}
\myS := \frac{1}{p} \sum_{\lambda_j\in I} (\hat{d}_j - \sigma_j^2) \toas 0.
\end{equation}
To borrow notation from \cite{erdHos2013averaging}, 
we conjecture further that, under the conditions of \cite{ledoit2011eigenvectors,ledoit2020analytical}, $\myS = O_\prec(p^{-d})$ uniformly in intervals $I$ for some $d>0$, where, $a_p = O_\prec(b_p)$ means that for any $\delta, D >0$, we have for sufficiently large $p$ that
\begin{equation} \label{eq:local}
\Pr\left[ |a_p| > p^{\delta} |b_p|  \right] < p^{-D}.
\end{equation} 
In other words, with polynomially high confidence, $a_p$ does not grow much more quickly than $b_p$, if at all. We further conjecture that
\begin{equation} \label{eq:microlocal}
\max_{i,j} \left| \bu_i'\bfr\bu_j - \delta_{ij} \hat{d}_j \right|  = O_\prec(p^{-d})
\end{equation}
for some $d > 0$.  This conjecture is motivated by an analogous result involving quadratic forms of Wigner-matrix eigenvectors \cite[Equation~1]{cipolloni2021eigenstate}.

Taking $c_i = \sigma_i^2/\hat{d}_i$, the test statistic in \eqref{eq:oZ0} becomes
\[
 \frac{1}{\sqrt{2\sum_i \sigma_i^4/\hat{d}_i^2}}\left(\tT^2 - \sum_{i=1}^p \frac{\sigma_i^2}{\hat{d}_i} \right).
\]
Using the fact that the above is asymptotically almost surely standard normal given $\bS$, and given the approximations of the last paragraph, we expect $Z$ to be asymptotically almost surely standard normal as well.  In particular, we expect $Z \gtrless \tau$ has an asymptotically constant false-alarm rate.


We present a simulation comparing a finite-sample distribution of $Z$ to a standard normal in Figure~\ref{fig:standard}.  Using 1000 Monte-Carlo values of $Z$ with $\bfr$ being the 200$\times$200 matrix $\bfr_4$ defined in the Section~\ref{sec:sim} and $n_1=n_2=200$, one can see that the finite-sample distribution of $Z$ approximates a standard normal, as predicted.




\begin{figure} 
\includegraphics[width=\columnwidth]{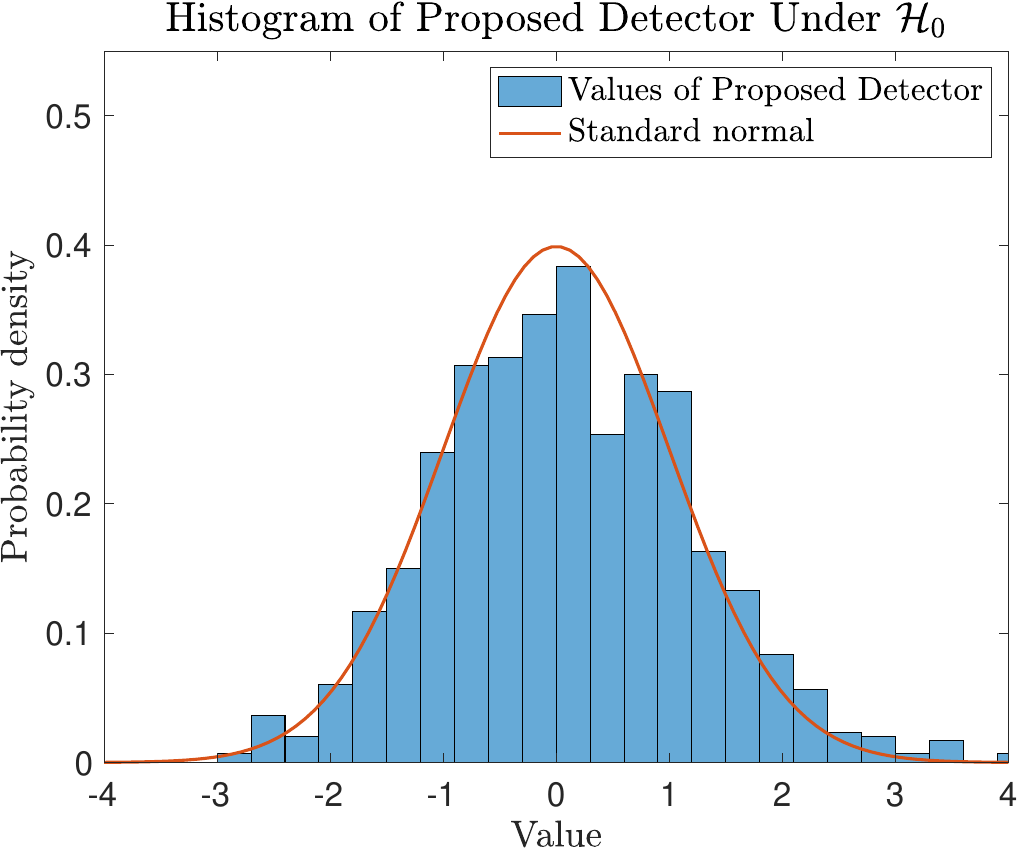}
\caption{A plot showing $Z$ from \eqref{eq:oZ} to be roughly standard normal.  Here, the population covariance matrix follows the model $\bfr_4$ in Section~\ref{sec:sim}.}
\label{fig:standard}
\end{figure}

\subsection{Maximizing Detection Power}


We now provide an intuitive argument that $\bfrhlw$ in \eqref{eq:ttilde0} optimizes conditional detection power given $\bS$ among shrinkage estimators $\bfrh$, for Gaussian data.  It can be shown using \eqref{eq:lwp-law} that $t = \tr(\bfr\bfrh^{-1})/p$ depends in the limit only on the limiting sample spectrum, so $\bfrh \leftarrow t\bfrh$ can be considered to be a shrinkage estimator, and we assume without loss of generality that $t$ converges to 1.

Consider the conditional false-alarm rate of the proposed detector given $\bS$: $\Pr\left[ \bd' \bfrh^{-1} \bd > \tau \mid \bS \right]$, where $\mathbf{d} = \bfr^{1/2}\bw$ and $\bw\sim \mathcal{N}(\mathbf{0},\mathbf{I})$. Using the result \cite[Theorem~5.1.4]{vershynin2018high} regarding concentration of Lipschitz functions on the sphere and the fact that $p\to\infty$, this probability can be approximated by
\begin{equation} \label{eq:cond-pfa}
\Pr\left[ \left\Vert\bw\right\Vert^2 \tr(\bfr^{1/2}\bfrh^{-1}\bfr^{1/2}) > \tau \mid \bS \right]. 
\end{equation}
Since $\bfrh$ depends only on $\bS$, the probability above is a function of $t$, $\tau$, and $p$ alone.  Since by assumption $t$ converges to 1, \eqref{eq:cond-pfa} is asymptotically dependent on $\tau$ and $p$ alone, and the same can be said of the corresponding unconditional probability.  Thus, we wish to maximize detection power for each $\tau$ and large $p$.

Consider the conditional probability of detection given $\bS$:
\begin{equation} \label{eq:cond-pd}
\Pr\left[(\bmu+\bd)'\bfrh^{-1}(\bmu+\bd)> \tau \mid \bS \right] ,
\end{equation}
where $\bmu =\sqrt{n_1 n_2/(n_1+n_2)} ( \bmu_1 - \bmu_2)$.
Expanding \eqref{eq:cond-pd}, we get
\begin{equation} \label{eq:pd2}
\Pr\left[\mya + \myb'\mathbf{v} + \left\Vert\bw\right\Vert^2 > \tau \mid \bS\right],
\end{equation}
where $\mya = \bmu'\bfrh^{-1}\bmu$,
$\myb = 2\bmu'\bfrh^{-1}\bfr^{1/2}\left\Vert\bw\right\Vert,$
and $\mathbf{v}=\bw/\left\Vert\bw\right\Vert$.  If we further condition on $\left \Vert \bw  \right \Vert$, the argument of the probability in \eqref{eq:pd2} depends only on the randomness of $\mathbf{v}$, and the probability itself depends only on $\alpha$ and $\left\Vert \myb \right\Vert$.  More precisely, the conditional probability of detection given $\bS$ and $\left \Vert \bw  \right \Vert$ is maximized when $\alpha/\left\Vert \myb\right\Vert$ is maximized.


For any $\left \Vert \bw  \right \Vert$ and $\bS$, maximizing $\mya/\left\Vert \myb\right\Vert$ is equivalent to
maximizing the \emph{signal-to-noise ratio}:
\begin{equation} \label{eq:snr0}
\frac{(\bmu'\bfrh^{-1}\bmu)^2}{\bmu'\bfrh^{-1}\bfr\bfrh^{-1}\bmu}.
\end{equation}
Using \cite[Theorem~5.1.4]{vershynin2018high} again, \eqref{eq:snr0} is well-approximated in probability for almost all $\bmu$ by
\begin{equation} \label{eq:snr}
    \frac{\tr(\bfrh^{-1})^2}{p\tr(\bfrh^{-1}\bfr\bfrh^{-1})},
\end{equation}
so that the ideal choice of $\bfrh$ subject to the constraint $t\approx 1$ is the Ledoit-Wolf estimator discussed in Section~\ref{sec:def}. This justifies our choice of detector.

An additional feature of equation~\eqref{eq:snr} is that it can be used to numerically find the optimum in the class of diagonal-loading estimators, which is a class of shrinkage estimators that includes \lalg{}.   This can be done by simply plugging $\bfrh = \bS + \lambda \mathbf{I}$ in \eqref{eq:snr} and iteratively solving for the optimal $\lambda > 0$ using, for example, \texttt{fminsearch} in Matlab.  Although this estimator is not technically the same as \lalg{}, we refer to it henceforth as \lalg{} since we can calculate it and its performance is an upper bound for \lalg{}'s.


\section{Simulations} \label{sec:sim}
 
  A key attribute of our proposed algorithm is its performance in simulation. In this section, we show performance that is consistent with average-case dominance over \bsalg{}, \cqalg{}, and \lalg{}, under the model that the difference in population means of the two samples is drawn uniformly from  the sphere.
 
In this Section, the dimension $p$ is $200$, $n_1$ and $n_2$ are taken to be 150, and we consider diagonal population covariance matrices $\bfr = \bfr_P$ parametrized by an integer $P \in \{0,2,4\}$.   The eigenvalues of $\bfr_P$ are defined for $1\le j \le 40$ by
 \[
 \left(\bfr_P\right)_{jj} = 10^{(41-j)P/40} + \epsilon_j, \qquad (1\le j \le 40),
 \]
 where $\epsilon_j$ is chosen i.i.d. uniformly at random from $[0,1]$.  On the other hand, for $j > 40$, we set  $(\bfr_P)_{jj} = 1$.
 Thus, largest eigenvalue of $\bfr_P$ is $10^P$, the next 39 decrease exponentially, and the remaining 160 are equal to unity.  This choice of spectrum mirrors the generalized spiked structure widely encountered in sensing applications, such as radar, and $P$ corresponds to the order of the $\bfr$'s condition number. 
 
 For each detector and each $P\in \{0,2,4\}$, we generate 100,000 detection scores as follows.  We generate a sub-Gaussian $p\times n_1$ and a $p \times n_2$ data matrix $\bX_1$ and $\bX_2$, both by coloring a matrix of i.i.d. uniform random variables with mean zero and variance one.  (Simulations appear to be similar for Gaussian data matrices.)  We  perturb the columns of one of these data matrices by a random vector uniformly drawn from a sphere of radius $1$ and treat the mean of the other as $\mathbf{0}$.  We then plug the data into the Proposed, \lalg{}, \bsalg{}, \cqalg{}, and Hotelling detectors.  
 The resulting 100,000 detection scores for each detector are used to generate the ROC curves in Figures~\ref{fig:R0}, \ref{fig:R2}, and \ref{fig:R4}, which correspond to $P=0,2,4$, respectively.  Figure~\ref{fig:R0} shows the case $P=0$, where \bsalg{}, \lalg{}, \cqalg{}, and the proposed method all perform similarly, and all dominate Hotelling.  In Figures~\ref{fig:R2} and \ref{fig:R4}, the proposed method dominates \bsalg{}, \lalg{}, \cqalg{}, and Hotelling, with the advantage being more pronounced for increasing values of $P$.  
 This is likely due to the fact that the number of large eigenvalues in these cases challenges both the spiked assumption of \lalg{} and the well-conditioned assumption of \bsalg{} and \cqalg{}, mentioned in Section~\ref{sec:background}.
 


\begin{figure}
\includegraphics[width=\columnwidth]{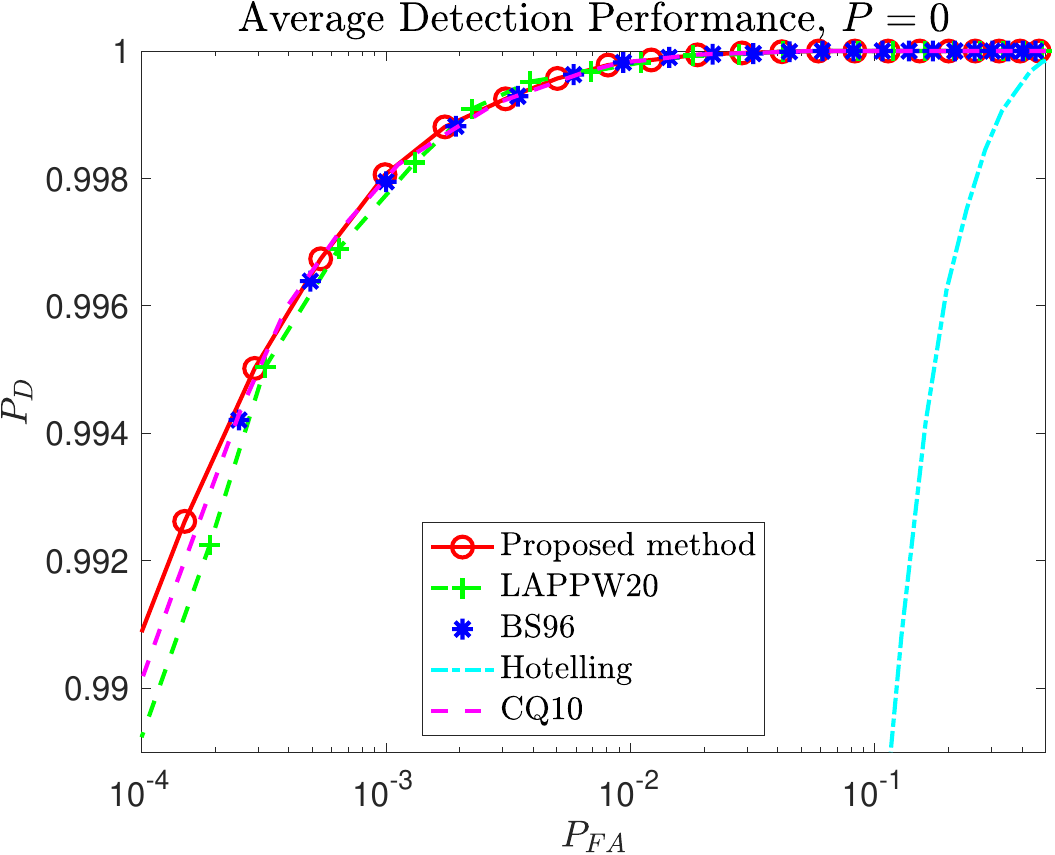}
\caption{For $\bfr= \bfr_0$ with $O(10^0)$ condition number, all methods perform similarly, except Hotelling.}
\label{fig:R0}
\end{figure}

\begin{figure}
\includegraphics[width=\columnwidth]{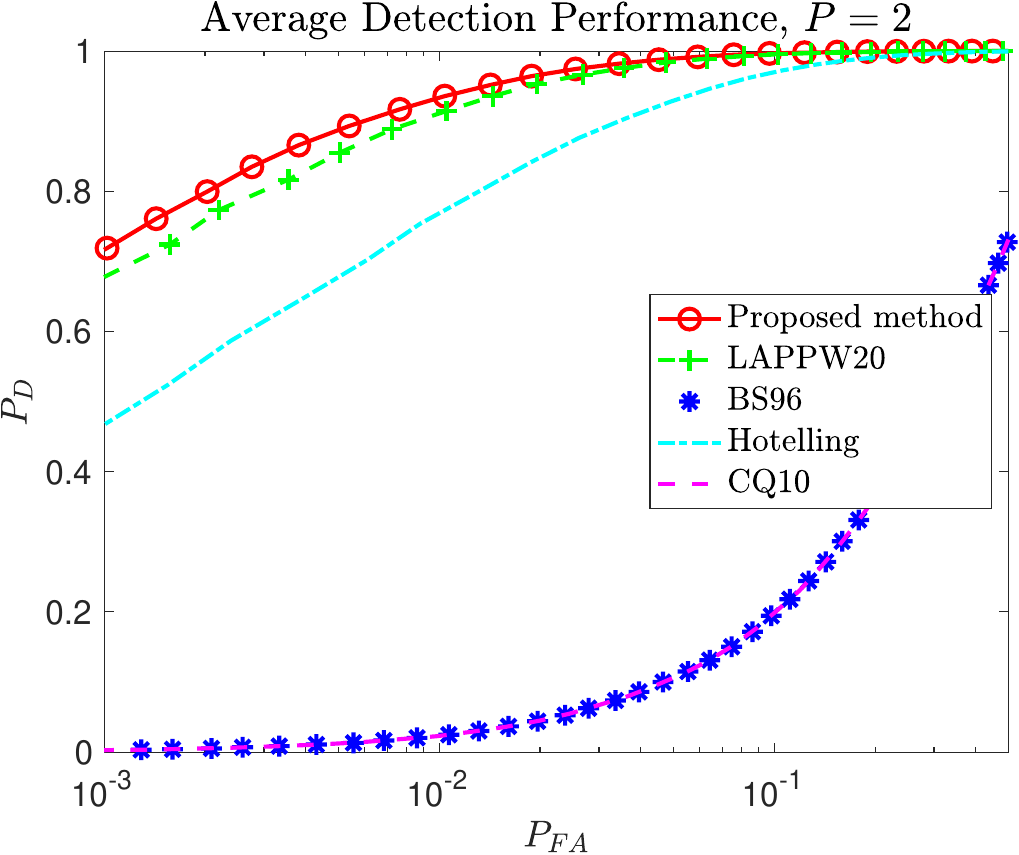}
\caption{For $\bfr= \bfr_2$ with $O(10^2)$ condition number, Proposed method performs similarly to \lalg{}, which outperforms Hotelling, which outperforms \bsalg{} and \cqalg{}, which lie on the chance line.}
\label{fig:R2}
\end{figure}


\begin{figure}
\includegraphics[width=\columnwidth]{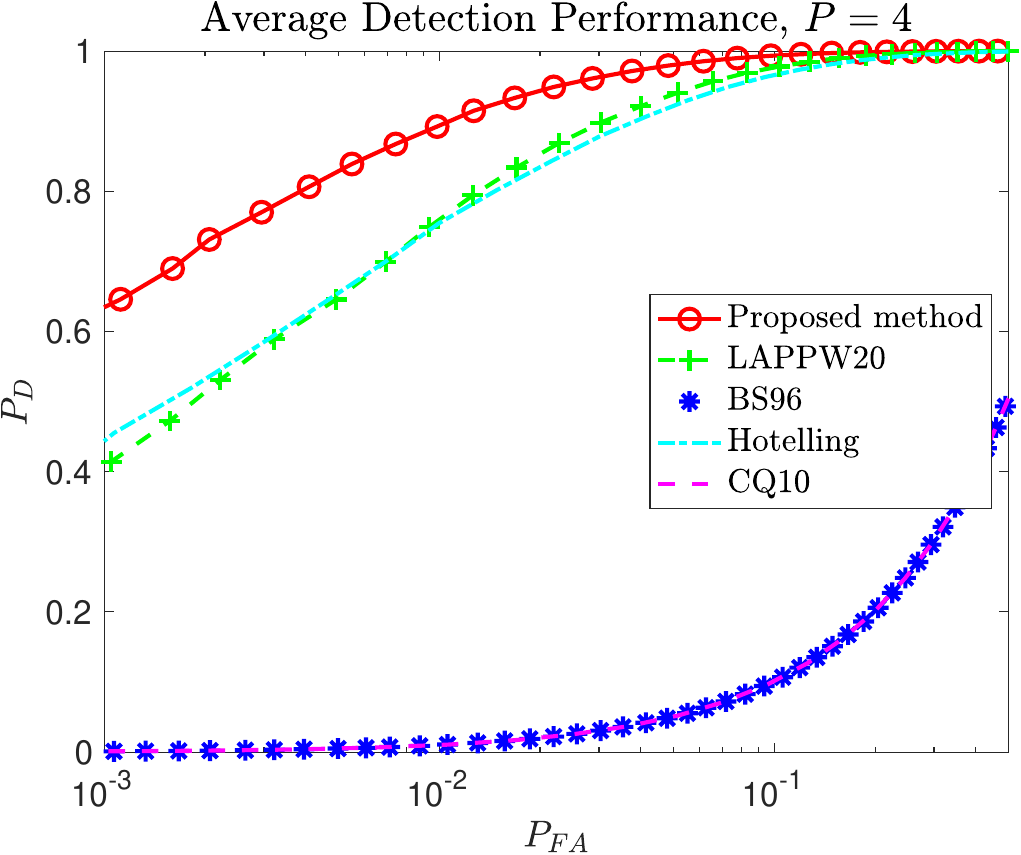}
\caption{For $\bfr= \bfr_4$ with $O(10^4)$ condition number, Proposed method outperforms \lalg{} and Hotelling, which outperform \bsalg{} and \cqalg{}, which lie on the chance line.}
\label{fig:R4}
\end{figure}

\section{Conclusion} \label{sec:conclusion}

In this paper, we have proposed an alternative to Hotelling's $T^2$ test with desirable theoretical and empirical properties in the high-dimensional setting.  In simulation, the proposed test appears to dominate the state-of-the-art alternatives \lalg{}, \cqalg{}, and \bsalg{}, likely due in part to the comparatively less forgiving assumptions they require. We provided an intuitive argument that our test has an asymptotically constant false-alarm rate and maximal detection power among shrinkage-based Hotelling-type tests under mild assumptions. We conjecture that a more rigorous foundation for this theory can be established by finding conditions under which \eqref{eq:microlocal} and strengthenings of \eqref{eq:lwp-law} hold, which we intend to investigate in future work.

\section*{Acknowledgments}

This work was supported by the United States Air Force Sensors Directorate, AFOSR grant 19RYCOR036, ARO grant W911NF-15-1-0479, and the DoD SMART Scholarship SEED grant. 
However, the views and opinions expressed in this article are those of the authors and do not necessarily reflect the official policy or position of any agency of the U.S. government.
Examples of analysis performed within this article are only examples. 
Assumptions made within the analysis are also not reflective of the position of any U.S. Government entity. The Public Affairs approval number of this document is AFRL-2022-2605.

\bibliographystyle{IEEEtran}
\bibliography{information-geometry-bib2}

\begin{thebibliography}{10}
\providecommand{\url}[1]{#1}
\csname url@samestyle\endcsname
\providecommand{\newblock}{\relax}
\providecommand{\bibinfo}[2]{#2}
\providecommand{\BIBentrySTDinterwordspacing}{\spaceskip=0pt\relax}
\providecommand{\BIBentryALTinterwordstretchfactor}{4}
\providecommand{\BIBentryALTinterwordspacing}{\spaceskip=\fontdimen2\font plus
\BIBentryALTinterwordstretchfactor\fontdimen3\font minus
  \fontdimen4\font\relax}
\providecommand{\BIBforeignlanguage}[2]{{%
\expandafter\ifx\csname l@#1\endcsname\relax
\typeout{** WARNING: IEEEtran.bst: No hyphenation pattern has been}%
\typeout{** loaded for the language `#1'. Using the pattern for}%
\typeout{** the default language instead.}%
\else
\language=\csname l@#1\endcsname
\fi
#2}}
\providecommand{\BIBdecl}{\relax}
\BIBdecl

\bibitem{anderson1963asymptotic}
T.~W. Anderson, ``Asymptotic theory for principal component analysis,''
  \emph{The Annals of Mathematical Statistics}, vol.~34, no.~1, pp. 122--148,
  1963.

\bibitem{paul2007asymptotics}
D.~Paul, ``Asymptotics of sample eigenstructure for a large dimensional spiked
  covariance model,'' \emph{Statistica Sinica}, pp. 1617--1642, 2007.

\bibitem{johnstone2001distribution}
I.~M. Johnstone, ``On the distribution of the largest eigenvalue in principal
  components analysis,'' \emph{Annals of Statistics}, pp. 295--327, 2001.

\bibitem{bai1996effect}
Z.~Bai and H.~Saranadasa, ``Effect of high dimension: {B}y an example of a two
  sample problem,'' \emph{Statistica Sinica}, pp. 311--329, 1996.

\bibitem{chen2010two}
S.~X. Chen and Y.-L. Qin, ``A two-sample test for high-dimensional data with
  applications to gene-set testing,'' \emph{The Annals of Statistics}, vol.~38,
  no.~2, pp. 808--835, 2010.

\bibitem{li2020adaptable}
H.~Li, A.~Aue, D.~Paul, J.~Peng, and P.~Wang, ``An adaptable generalization of
  {H}otelling's {T}$^2$ test in high dimension,'' \emph{The Annals of
  Statistics}, vol.~48, no.~3, pp. 1815--1847, 2020.

\bibitem{ledoit2020analytical}
O.~Ledoit and M.~Wolf, ``Analytical nonlinear shrinkage of large-dimensional
  covariance matrices,'' \emph{The Annals of Statistics}, vol.~48, no.~5, pp.
  3043--3065, 2020.

\bibitem{stein1975estimation}
C.~Stein, ``Estimation of a covariance matrix, {R}ietz lecture,'' in \emph{39th
  Annual Meeting IMS, Atlanta, GA, 1975}, 1975.

\bibitem{stein1986lectures}
------, ``Lectures on the theory of estimation of many parameters,''
  \emph{Journal of Soviet Mathematics}, vol.~34, no.~1, pp. 1373--1403, 1986.

\bibitem{robinson2021space}
B.~D. Robinson, R.~Malinas, and A.~O. Hero, ``Space-time adaptive detection at
  low sample support,'' \emph{IEEE Transactions on Signal Processing}, vol.~69,
  pp. 2939--2954, 2021.

\bibitem{namdari2021high}
J.~Namdari, D.~Paul, and L.~Wang, ``High-dimensional linear models: {A} random
  matrix perspective,'' \emph{Sankhya A}, vol.~83, no.~2, pp. 645--695, 2021.

\bibitem{donoho2018optimal}
D.~L. Donoho, M.~Gavish, and I.~M. Johnstone, ``Optimal shrinkage of
  eigenvalues in the spiked covariance model,'' \emph{Annals of Statistics},
  vol.~46, no.~4, p. 1742, 2018.

\bibitem{ledoit2011eigenvectors}
O.~Ledoit and S.~P{\'e}ch{\'e}, ``Eigenvectors of some large sample covariance
  matrix ensembles,'' \emph{Probability Theory and Related Fields}, vol. 151,
  no. 1-2, pp. 233--264, 2011.

\bibitem{ledoit2012nonlinear}
O.~Ledoit and M.~Wolf, ``Nonlinear shrinkage estimation of large-dimensional
  covariance matrices,'' \emph{The Annals of Statistics}, vol.~40, no.~2, pp.
  1024--1060, 2012.

\bibitem{ledoit2018optimal}
------, ``Optimal estimation of a large-dimensional covariance matrix under
  {S}tein's loss,'' \emph{Bernoulli}, vol.~24, no.~4B, pp. 3791--3832, 2018.

\bibitem{benaych2016lectures}
F.~Benaych-Georges and A.~Knowles, ``Lectures on the local semicircle law for
  {W}igner matrices,'' \emph{arXiv preprint arXiv:1601.04055}, 2016.

\bibitem{erdHos2013averaging}
L.~Erd{\H{o}}s, A.~Knowles, and H.-T. Yau, ``Averaging fluctuations in
  resolvents of random band matrices,'' in \emph{Annales Henri Poincar{\'e}},
  vol.~14, no.~8.\hskip 1em plus 0.5em minus 0.4em\relax Springer, 2013, pp.
  1837--1926.

\bibitem{cipolloni2021eigenstate}
G.~Cipolloni, L.~Erd{\H{o}}s, and D.~Schr{\"o}der, ``Eigenstate thermalization
  hypothesis for {W}igner matrices,'' \emph{Communications in Mathematical
  Physics}, vol. 388, no.~2, pp. 1005--1048, 2021.

\bibitem{vershynin2018high}
R.~Vershynin, \emph{High-Dimensional Probability: An Introduction with
  Applications in Data Science}.\hskip 1em plus 0.5em minus 0.4em\relax
  Cambridge University Press, 2018, vol.~47.

\end{thebibliography}

\end{document}